\documentclass[11pt]{amsart}
\usepackage{graphicx}
\usepackage{amsfonts, amsmath, amsfonts, amssymb}
\newtheorem{thm}{\rm THEOREM}[section]

\newtheorem{lem}[thm]{\rm LEMMA}

\theoremstyle{definition}
\newtheorem{defn}[thm]{\it Definition}
\theoremstyle{remark}
\newtheorem{rem}[thm]{\it Remark}
\numberwithin{equation}{section}

\numberwithin{equation}{section}

\numberwithin{equation}{section} \theoremstyle{cond}

 \numberwithin{equation}{section}
 
\newenvironment{prf}{ \noindent{\bf Proof}}{\\ \hspace*{\fill}$\Box$ \par  }

\begin{document}


\title[]{ Intrinsic Differential Geometry and the Existence of
Quasimeromorphic Mappings}
\author{Emil Saucan}%
\address{\it Mathematics Department, Technion, Technion City, Haifa 32000 }%
\email{semil@tx.technion.ac.il}%


\date{\today}%


\begin{abstract}
We give a new proof of the existence of nontrivial quasimeromorphic
mappings on a smooth Riemannian manifold, using solely the intrinsic
geometry of the manifold.
\end{abstract}

\maketitle


\section{Introduction and Background}

The existence of quasimeromorphic ($qm$) mappings on
$\mathcal{C}^\infty$-Riemannian manifolds without boundary is due to
Peltonen \cite{pe}, and represents a generalization of previous
results of Tukia \cite{tu} and Martio-Srebro \cite{ms}. In \cite{s1}
we have extended Peltonen's result to include manifolds with
boundary and of lower differentiability class. A further
generalization to certain classes of orbifolds was given in
\cite{s2}, \cite{s3}.

The essential ingredient in all the results above is construction of
a {\it thick} (or {\it fat}) ``chessboard triangulation'' (i.e. such
that two given $n$-simplices having a $(n-1)$-dimensional face in
common will have opposite orientations), each of its simplices being
then quasiconformally mapped on the unit sphere $\mathbb{S}^n$ using
the classical {\it Alexander method} \cite{al}.

Recall that thick triangulations are defined as follows:

\begin{defn} Let $\tau \subset \mathbb{R}^n$ ; $0 \leq k \leq n$ be a $k$-dimensional simplex.
The {\it thickness}  $\varphi$ of $\tau$ is defined as being:
\begin{equation}
\varphi = \varphi(\tau) = \hspace{-0.3cm}\inf_{\hspace{0.4cm}\sigma
< \tau
\raisebox{-0.25cm}{\hspace{-0.9cm}\mbox{\scriptsize$dim\,\sigma=j$}}}\!\!\frac{Vol_j(\sigma)}{diam^{j}\,\sigma}\;.
\end{equation}
The infimum is taken over all the faces of $\tau$, $\sigma < \tau$,
and $Vol_{j}(\sigma)$ and $diam\,\sigma$ stand for the Euclidian
$j$-volume and the diameter of $\sigma$ respectively. (If
$dim\,\sigma = 0$, then $Vol_{j}(\sigma) = 1$, by convention.)

 A simplex $\tau$ is $\varphi_0${\it-thick}, for some $\varphi_0 > 0$,
if $\varphi(\tau) \geq \varphi_0$. A triangulation (of a submanifold
of $\mathbb{R}^n$) $\mathcal{T} = \{ \sigma_i \}_{i\in \bf I}$ is
$\varphi_0${\it-thick} if all its simplices are $\varphi_0$-thick. A
triangulation $\mathcal{T} = \{ \sigma_i \}_{i\in \bf I }$ is {\it
thick} if there exists $\varphi_0 \geq 0$ such that all its
simplices are $\varphi_0${\it-thick}.
\end{defn}

The definition above is the one introduced in \cite{cms}. For some
different, yet equivalent definitions of thickness, see \cite{ca1},
\cite{ca2}, \cite{mun}, \cite{pe}, \cite{tu}.

Note that in our generalizations \cite{s1}, \cite{s2} we have used
Peltone's result, in conjunction with methods of Munkres \cite{mun}
and Cheeger et al. \cite{cms} to obtain the desired thick
triangulation.

The method of proof employed in \cite{pe} is based upon {\it
extrinsic} Differential Geometric considerations. %
More precisely, the idea of the proof is as follows: Start by
isometrically embedding the $n$-dimensional, complete, Riemannian
manifold $M^n$ into $\mathbb{R}^\nu$, for some large enough $\nu$.
(The existence of such an embedding dimension ``$\nu$'' is
guaranteed by Nash's Embedding Theorem \cite{na}.) Then one
constructs an exhaustion of $M^n$ by a sequence of compact manifolds
$\{M_i\}_{i \in \mathbb{N}}$\,.

To control the size of these compact manifolds and that of the
``pasting zones'' between them (as well as the density of the set of
vertices of the triangulation to be constructed), one makes appeal
to two geometric features, namely the {\em osculatory} (or {\em
tubular}) {\em radius} and the {\em connectivity radius}, who are
defined as follows:

\begin{defn}
\begin{enumerate}
\item $\mathbb{S}^{\nu-1}(x,\rho)$ is an {\rm osculatory sphere} at $x \in M^n$ iff:
\begin{enumerate}
\item $\mathbb{S}^{\nu-1}(x,\rho)$ is tangent at x;
\\ and
\item $\mathbb{B}^n(x,\rho) \cap M^n = \emptyset$.
\end{enumerate}
\item Let $X \subset M^n$. The number $\omega = \omega_X = \sup\{\rho > 0\,|\, \mathbb{S}^{\nu-1}(x,\rho) \; {\rm osculatory} \\{\rm at\; any}\; x \in
X\}$ is called the {\em maximal osculatory} ({\em tubular}) {\em
radius} at $X$.
\end{enumerate}
\end{defn}
\hspace*{-0.4cm}where tangentiality generalizes in a straightforward
manner the classical notion defined for surfaces in $\mathbb{R}^3$:

\begin{defn}
$\mathbb{S}^{\nu-1}(x,r)$ is {\em tangent} to $M^n$ at $x\in M^n$
iff there exists $\mathbb{S}^n(x,r) \subset
\mathbb{S}^{\nu-1}(x,r)$, such that $T_x(\mathbb{S}^n(x,r)) \equiv
T_x(M^n)$.
\end{defn}
\hspace*{-0.4cm}(Here $\mathbb{B}^\nu(x,r) = \{y \in
\mathbb{R}^\nu\,|\, d_{eucl} < r\}$; $\mathbb{S}^{\nu-1}(x,r)=
\partial\mathbb{B}^\nu(x,r)$.)

Note that there exists an osculatory sphere at any point of $M^n$
(see, e.g. \cite{pe}\,).

\begin{defn} Let $U \subset M^n, U \neq \emptyset$, be a relatively compact set, and let $T = \bigcup_{x \in
\bar{U}}\sigma(x,\omega_U)$. The number $\kappa_U =
\max\{r\,|\,\sigma^n(x,r)\;  {\rm  is\; connected \; for \; all}\; s
\leq \omega_U,\, x \in \bar{T}\}$, is called the {\em maximal
connectivity radius} at U.
\end{defn}
\hspace*{-0.4cm}(Here and $\sigma^n(x,r) = M^n \cap
\mathbb{B}^\nu(x,r)$\,.)

These geometric features help us assure that the manifold does not
``turn on itself too fast'', piercing a simplex of the (future)
triangulation. Moreover, they are interrelated through the following
inequality (see \cite{pe}, Lemma 3.1):

\[\omega_U \leq \frac{\sqrt{3}}{3}\kappa_U\,.\] \label{ec:1}

It follows, therefore, that to obtain a vertex set of the required
 density, one can employ estimates that are
solely functions of $\omega_U$.

Obviously, this construction is basically extrinsic, since it
essentially uses the Nash embedding and because the geometric
features that control the density of the vertices and the thickness
of the simples are also extrinsic (see above). Moreover, computing
the osculatory and connectivity radii is very difficult. Even
computing the principal curvatures of the Nash embedding by solving
the specific Gauss Equation is highly problematic. (See \cite{s3}
for a discussion of these aspects and also \cite{saz} for their
applicative side.)


We have recently given in \cite{SK} a simpler proof of the existence
of thick triangulations on manifolds (and hence of $qm$-mappings),
where by ``simpler'' we mean that it mainly uses tools of Elementary
Differential Topology. However, this proof still requires the
embedding of $M^n$ into some $\mathbb{R}^N$, for $N$ large enough,
 hence it is still partially extrinsic
in nature. More important, most of the geometric information
regarding the manifold is (evidently) lost or hard to retrieve when
using the Differential Topology approach. However, in many cases the
manifold comes not merely endowed with a Riemannian metric, but also
with some more concrete information on its geometry, usually in the
form of bounds for curvatures, volume and diameter. It is therefore
useful to have a construction that uses this geometric data. It is
the goal of this paper to produce precisely such a construction,
which we present in the next section. Finally, for the sake of
completeness, in the last section we remind the reader how a
quasimeromorphic mapping is obtained once a thick triangulation is
constructed.



\section{The construction}

As in Peltonen's construction, the idea of the proof 
is to use the basic fact that $M^n$ is $\sigma$-compact, i.e. it
admits an exhaustion by compact submanifolds $\{M_i\}_i$ (see, e.g.
\cite{spi}). This is a standard fact for metrizable manifolds.
However, it is conceivable that the ``cutting surfaces'' $N_{ij}$\,,
$\bigcup_{\scriptscriptstyle{j=1,...k_i}}\hspace{-0.2cm}N_{ij} =
\partial M_j$\,, are merely $\mathcal{C}^0$, so even the existence of a
triangulation for these hypersurfaces is not always assured, 
hence a fortiori that of smooth triangulations. (See. e.g. \cite{th}
for a brief review of the results regarding the existence of
triangulations).

To show that one can obtain (by ``cutting along'') smooth
hypersurfaces, we briefly review the main idea of the proof of the
$\sigma$-compactness of $M^n$ (for the full details, see, for
example \cite{spi}): Starting from an arbitrary base point $x_0 \in
M^n$, one considers the interval $I = I(x_0) = \{r > 0\,|\,
\beta^n(x_0,r)\; {\rm is \; compact}\}$; $\beta^n(x,r) =
exp_x\big(\mathbb{B}^n(0,r)\big)$, where $exp_x$ denotes the
exponential map: $exp_x:T_x(M^n) \rightarrow M^n$, and where
$\mathbb{B}^n(0,r) \subset T_x\big(M^n\big)$, $\mathbb{B}^n(0,r) =
\{y \in \mathbb{R}^n\,|\,d_{eucl}(y,0) < r\}$. If $I = \mathbb{R}$,
then $M^n = \bigcup_{\scriptscriptstyle
1}^{\scriptscriptstyle\infty}{\beta^n(x,i)}$, thence
$\sigma$-compact. If $I \neq \mathbb{R}$, one constructs the
compacts sets $M_i$, $M_0 = \{x_0\}$, $M_{i+1} =
\bigcup_{\scriptscriptstyle y \in M_i}\beta^n(y,r(y))$, where $r(y)
= \frac{1}{2}\sup\{r \in I(y)\}$. Then it can be shown that $M^n =
\bigcup_{\scriptscriptstyle n \geq 0}M_i$, i.e. $M^n$ is
$\sigma$-compact.

The smoothness of the surfaces $N_{ij}$ now follows from Wolter's
result \cite{wo} regarding the $2$-differentiability of the cut
locus of the exponential map.

We shall construct 
thick triangulations of $M_i$ and $N_{ij}$ of thickness $\varphi_1 =
\varphi_1(n)$ and $\varphi_2 = \varphi_2(n-1)$, respectively. We can
then apply repeatedly the ``mashing'' technique developed in
\cite{s2}, for collars of $N_{ij}$ in $M_i$ and $M_{i + 1}$, $j \geq
0$, rendering a triangulation of $M^n$, of uniform thickness
 $\varphi = \varphi(n)$ (see \cite{s2}, \cite{cms}).

Up to this point, our construction is practically identical to that
we used in \cite{SK}. However, to produce the fat triangulations of
$M_i$ and $N_{ij}$, we shall employ, as stated before, methods of
Intrinsic Differential Geometry, instead of(rather than) the ones of
Differential Topology we applied in \cite{SK}.

We start by noting that the manifolds $M_i$, and $N_{ij}$\,, $i,j
\in \mathbb{N}$ are compact, hence the have bounded {\it sectional
curvatures} (see, e.g. \cite{Be}) and diameters. Let $k_i$,$k_{ij}$
and $K_i,K_{ij}$ denote the lower bound, respective the upper bound,
for the sectional curvatures, and let $D_i$,$D_{ij}$ denote the
upper bound of the diameter of $M_i$ and $N_{ij}$, respectively.

Therefore, for each of these manifolds, we can make avail of a
triangulation method that,  according to \cite{Be}, was developed,
yet not published, by Karcher, but which, to the best of our
knowledge, appeared for the first time in \cite{GP}. (The same
method was applied by Weinstein \cite{We}, to obtain a similar
result in even dimension.)

The idea is to use so called {\it efficient packings}:

\begin{defn}
Let $p_1,\ldots,p_{n_0}$ be points  $\in M^n$,  satisfying  the
following conditions:
\begin{enumerate}
\item The set $\{p_1,\ldots,p_{n_0}\}$ is an $\varepsilon$-net on $M^n$, i.e. the
balls $\beta^n(p_k,\varepsilon)$, $k=1,\ldots,n_0$ cover $M^n$;
\item The balls (in the intrinsic metric of $M^n$) $\beta^n(p_k,\varepsilon/2)$ are pairwise
disjoint.
\end{enumerate}
Then the set $\{p_1,\ldots,p_{n_0}\}$ is called a {\it minimal
$\varepsilon$-net} and the packing with the balls
$\beta^n(p_k,\varepsilon/2)$, $k=1,\ldots,n_0$, is called an {\it
efficient packing}. The set $\{(k,l)\,|\,k,l = 1,\ldots,n_0\; {\rm
and}\; \beta^n(p_k,\varepsilon) \cap \beta^n(p_l,\varepsilon) \neq
\emptyset\}$ is called the {\it intersection pattern} of the minimal
$\varepsilon$-net (of the efficient packing).
\end{defn}

Efficient packings have the following important properties, which we
list below (for proofs see \cite{GP}):

\begin{lem}
There exists $n_1 = n_1(n,k_i,D_i)$, such that if
$\{p_1,\ldots,p_{n_0}\}$ is an $\varepsilon$-net on $M^n$, then $n_0
\leq n_1$.
\end{lem}

\begin{lem}
There exists $n_2 = n_2(n,k_i,D_i)$, such that for any $x \in M^n$,
$\left|\{k \,|\, k = 1,\ldots,n_0\; {\rm and}\;
\beta^n(x,\varepsilon) \cap \beta^n(p_k,\varepsilon) \neq
\emptyset\}\right| \leq n_2$, for any minimal $\varepsilon$-net
$\{p_1,\ldots,p_{n_0}\}$.
\end{lem}

\begin{lem}
Let $M^n, \mathfrak{M}^n$, be manifolds having the same bounds $k_i$
and $D_i$ (see above) and let $\{p_1,\ldots,p_{n_0}\}$ and
$\{q_1,\ldots,q_{n_0}\}$ be minimal $\varepsilon$-nets with the same
intersection pattern, on $M^n$, $\mathfrak{M}^n$, respectively. Then
there exists a constant $n_3 = n_3(n,k_i,D_i,K_i)$, such that if
$d(q_i,q_j) < K_i\cdot\varepsilon$, then $d(q_i,q_j) <
n_3\cdot\varepsilon$.
\end{lem}

Such an efficient packing is always possible on a closed, connected
Riemannian manifold and, by using the properties above, one can
construct a simplicial complex
 having as vertices the centers of the balls
 $\beta^n(p_k,\varepsilon)$. (Edges are connecting the centers of
 adjacent balls; further edges being added to ensure the cell
 complex obtained is triangulated to obtain a simplicial complex.)

\begin{rem}
Let $M^n, \mathfrak{M}^n$, be manifolds having the same bounds $k,
D$ and $v$, where $v$ denotes the lower bound for volume. There
exists an $\varepsilon = \varepsilon(k,D,v)
> 0$, such that any two minimal $\varepsilon$-nets on $M^n,
\mathfrak{M}^n$ with the same intersection pattern, are
homeomorphic. Moreover, given $k,D$ and $v$ as above, the number of
such homotopy types is finite (see \cite{GP}).
\end{rem}

One can ensure that the triangulation will be convex and that its
simplices are convex, by choosing $\varepsilon = {\rm
ConvRad}(M^n)$, where the {\it convexity radius} ${\rm
ConvRad}(M^n)$ is defined as follows:

\begin{defn}
Let $M^n$ be a Riemannian manifold. The {\it convexity radius} of
$M^n$ is defined as $\inf\{r>0\,|\, \beta^n(x,r) \; {\rm is\;
convex},\; {\rm for\; all\;} x \in M^n\}$.
\end{defn}

This follows from the fact that $\beta^n\left(x,{\rm
ConvRad}(M^n)\right) \subset \beta^n\left(x,{\rm
InjRad}(M^n)\right)\,$, (since ${\rm ConvRad}(M^n) = \frac{1}{2}{\rm
InjRad}(M^n)$ -- see, e.g.  \cite{Be}). Here ${\rm InjRad}(M^n)$
denotes the {\it injectivity radius}:

\begin{defn} Let $M^n$  be a Riemannian manifold. The {\em
injectivity radius} of $M^n$ is defined as: ${\rm InjRad}(M^n) =
\inf{x \in M^n}\,|\,{Inj(x)}$, where $Inj(x) = \sup\{{r\,|\,
exp_{x}|_{\mathbb{B}^n(x,r)}\; {\rm is\; a \; diffeomorphism}}\}$.
\end{defn}

Note that by a classical result of Cheeger \cite{Ch}, there is a
universal positive lower bound for ${\rm InjRad}(M_i)$ in terms of
$k_i, D_i$ and $v_i$, where $v_i$ is the lower bound for the volume
of $M_i$. It is precisely this result (and similar ones -- see also
the discussion below) that make the triangulation exposed above a
simple and practical one, at least in many cases.

The same method of triangulation can be applied to the manifolds
$N_{ij}$. The triangulations thus obtained can be ``thickened'' by
applying the techniques of \cite{cms} or \cite{s4}. Then one can
``mash'' and ``thicken'' the triangulations of $M_i$ and $N_{ij}$,
to obtain using the method of \cite{s1}. Applying this process
inductively for all the elements of the exhaustion $\{M_i\}_{i \in
\mathbb{N}}$, one obtains a uniformly thick triangulation of $M^n$,
thus concluding the announced alternative proof of Peltonen's
result:

\begin{thm}  \label{thm:pelt}
Let $M^n$, $n \geq 2$, be complete, connected, $\mathcal{C}^\infty$
Riemannian manifold. Then $M^n$ admits a (uniformly) thick
triangulation.
\end{thm}

\begin{rem}
As already noted in the introduction, the result above can be
extended to manifolds with boundary, of low differentiability class
(see \cite{s1}) and to certain types of orbifolds (see \cite{s2},
\cite{s3}).
\end{rem}

\begin{rem} \label{rem:abb}
Instead of using the method of estimating convexity radii for the
manifolds $int M_i$ and their boundary components $N_{ij}$, we could
have used the estimates for the convexity radii of $N_{ij}$ using
the methods of \cite{ABB}. However, this would have provided more
difficult. In addition, the hierarchical approach adopted here is
the classical one of \cite{ca1}.
\end{rem}

\begin{rem} \label{rem:ric}
The same basic method of triangulation as employed herein may be
applied by considering bounds on the {\it Ricci curvature} of the
manifolds $M_i$ and $N_{ij}$\,, $i,j \in \mathbb{N}$. This
relaxation allows us to apply this technique to manifolds for which
less geometric control is possible.
\end{rem}

We conclude this section by reviewing the advantages and
disadvantages of the triangulation method introduced above, as
compared to that of \cite{pe}. As we have already noted above, an
immediate advantage stems from the fact that, by using solely the
intrinsic geometry of the manifold, this approach does not
necessitate the cumbersome Nash embedding technique, that results in
a quasi-impossible computation of the curvatures required in
Peltone's construction. But, perhaps, the main advantage resides in
the fact that {\it universal} lower bounds can be computed for the
injectivity radius (hence also for the convexity radius). Besides
Cheeger's classical result mentoned above, many other such theorems
for compact manifolds exist -- see \cite{Be} for
a plethora of relevant theorems. 
In addition, as we mentioned in Remark \ref{rem:abb}, a lower bound
for injectivity radius can be determined also for manifolds with
boundary. Moreover, as noted in Remark \ref{rem:ric}, such bounds
can be attained in terms of the Ricci curvature, further extending
the class of the manifolds for which our method may be easily
applied. It should also be noted that, by a theorem of Maeda
\cite{Ma}, for certain types of noncompact manifolds a universal
lower bound also exist, more precisely for (noncompact) manifolds
with sectional curvatures $K$ satisfying $0 < K \leq K_0$, the
following inequality holds: ${\rm InjRad}(M_i) \geq \pi/\sqrt{K_0}$.
Moreover, if $n=2$ and if $M^n$ is homeomorphic to $\mathbb{R}^2$,
then the same lower bound is achieved even under the weaker
assumption that $0 \leq K \leq K_0$.

The main disadvantage, as compared to \cite{pe}, of the approach
adopted herein, resides in the lack of control of the curvatures of
the ``cutting'' surfaces $N_{ij}$. In consequence, possible drastic
changes in sectional curvatures may occur, thence in injectivity
radii and, implicitly, in the size of the simplices, when passing
from $M_i$ to $N_{ij}$. (As a typical case for this kind of
behavior, consider a ``crumpled'' closed ball in $\mathbb{R}^3$.
Then its interior has the trivial Euclidean geometry of the ambient
space, whence sectional curvatures  $\equiv 0$, while the sectional
(i.e. Gaussian) curvature of the boundary may attain arbitrarily
large values of $|K|$.)

However, it is possible to smoothen the Riemannian metric of $M^n$,
to obtain a metric having a sectional curvatures bound, while
remaining arbitrarily close to the original metric (see, e.g.
\cite{CG}, \cite{Y}). (We should note here, in conjunction with
Remark \ref{rem:ric}, that results regarding the smoothing of
sectional curvature, under weaker Ricci curvatures bounds, also
exist (see, e.g. \cite{DWY}, \cite{an}, \cite{PWY}).


\section{The existence of quasimeromorphic mappings}

We begin this section by reminding the reader the definition of
quasimeromorphic mappings:

\begin{defn}
Let $M^n, N^n$ be oriented, Riemannian $n$-manifolds.
\begin{enumerate}
\item $f:M^n \rightarrow N^n$ is called {\it quasiregular} ($qr$) iff
\begin{enumerate}
\item $f$ is locally Lipschitz (and thus differentiable a.e.);
\\ \hspace*{-0.8cm} and
\item \(0 < |f'(x)|^n \leq KJ_f(x)\), for any \(x \in M^n\);
\end{enumerate}
\;where $f'(x)$ denotes the formal derivative of $f$ at $x$,
$|f'(x)| = \sup
\raisebox{-0.25cm}{\mbox{\hspace{-0.75cm}\tiny$|h|=1$}}|f'(x)h|$,
and where $J_f(x) = detf'(x)$;
\item {\it quasimeromorphic} ($qm$) iff $N^n = \mathbb{S}^n$,
where $\mathbb{S}^n$ is usually identified with
$\widehat{\mathbb{R}^n} = \mathbb{R}^n \cup \{\infty\}$ endowed with the {\it spherical metric}. 
\end{enumerate}
The smallest number $K$ that satisfies condition (b) above is called
the {\it outer dilatation} of \nolinebreak[4]$f$.
\end{defn}

Before proceeding further, we need the following technical lemma
(for its proof, see \cite{ms}, \cite{pe}):

\begin{lem}
Let $\mathcal{T}$ be a fat triangulation of $M^n \subset
\mathbb{R}^N$, and let $\tau,\sigma \in
\nolinebreak[4]\mathcal{T},\; \tau = (p_0,\dots,p_n), \,\sigma =
(q_0,\dots,q_n)$; and denote $|\tau| = \tau \cup int\,\tau$. Then
there exists a orientation-preserving homeomorphism $h = h_{\tau}:
|\tau| \rightarrow \widehat{\mathbb{R}^{n}}$ such that:
\begin{enumerate}
\item $h(|\tau|) = |\sigma|$, \,if\; $det(p_0,\dots,p_n) > 0$
\\ and
\\ $h(|\tau|) = \widehat{\mathbb{R}^{n}} \setminus| \sigma|$, \,if\; $det(p_0,\dots,p_n) < 0$.
\item $h(p_i) = q_i, \; i=0,\ldots,n.$
\item $h|_{\partial|\sigma|}$ is a $PL$ homeomorphism.
\item $h|_{int|\sigma|}$ is quasiconformal.
\end{enumerate}
\end{lem}

We can now prove the existence of $qm$-mappings on open Riemannian
manifolds:

\begin{thm}
Let $M^n, n \geq 2$, be a connected, complete, oriented $C^\infty$
Riemannian manifold. 
Then there exists a non-constant quasimeromorphic mapping $f:M^n
\rightarrow \widehat{\mathbb{R}^n}$.
\end{thm}

\begin{prf} Let $\mathcal{T}$ be the thick triangulation provided by Theorem \ref{thm:pelt}. 
Furthermore, by performing a
barycentric type subdivision before starting the fattening process
of the triangulation given by Theorem \ref{thm:pelt}, ensure that
all the simplices of the triangulation satisfy the condition that
every $(n-2)$-face is be incident to an even number of
$n$-simplices. Let $f:M ^n\rightarrow \widehat{\mathbb{R}^n}$ be
defined by: $f|_{|\sigma|} = h_{\sigma}$, where $h$ is a
homeomorphism constructed in the lemma above. Then $f$ is a local
homeomorphism on the $(n-1)$-skeleton of $\mathcal{T}$ too, while
its branching set $B_f$ is the $(n-2)$-skeleton of $\mathcal{T}$. By
its construction $f$ is quasiregular. Moreover, given the uniform
fatness of the triangulation $\mathcal{T}$, the dilatation of $f$
depends only on the dimension \nolinebreak[4]$n$.
\end{prf}

\begin{rem}
Again, this result may be extended to include manifolds with
boundary, of low differentiability class (see \cite{s1}) and to
certain types of orbifolds (see \cite{s2}, \cite{s3}).
\end{rem}


\subsection*{Acknowledgment}

The author wishes to express his gratitude to Professor Shahar
Mendelson for his warm support and for his stimulating questions.
He would also like to thank 
Professor Klaus-Dieter Semmler for bringing to his attention
Wolter's work and to Professor Meir Katchalski, who set things into
motion.




\end{document}